\documentclass[greek,english]{article}
\usepackage{babel}
\usepackage{amssymb}
\usepackage{amsmath,amsthm}
\usepackage{amsxtra}
\usepackage[mathscr]{eucal}
\usepackage{graphics}
\usepackage{epsfig}

\renewcommand{\le}{\leqslant}
\renewcommand{\ge}{\geqslant}

\makeatletter                   
\renewcommand\@biblabel[1]{#1.} 
\makeatother

\renewcommand\mid{:}

\newtheorem{theorem}{Theorem}[section]
\newtheorem{lemma}[theorem]{Lemma}

\newtheorem{proposition}[theorem]{Proposition}

\theoremstyle{definition}

\newtheorem{conjecture}[theorem]{Conjecture}
\newtheorem{example}{Example}

\def\N{\mathbb{N}}

\def\primeish{biased}

\def\cB{{\mathcal B}}
\def\cR{{\mathcal R}}
\def\cS{{\mathcal S}}
\def\cT{{\mathcal T}}
\def\cU{{\mathcal U}}

\def\ord{\operatorname{ord}}

\begin{document}
\title{Primes Generated by Recurrence
Sequences}
\author{Graham Everest, Shaun Stevens, Duncan Tamsett, and Tom Ward}

\maketitle

\addtocounter{section}{1}
\noindent{\bf\thesection\ MERSENNE NUMBERS AND PRIMITIVE PRIME
DIVISORS.}\setcounter{theorem}{0} A
notorious problem from elementary number theory is the ``Mersenne
Prime Conjecture.'' This asserts that the \emph{Mersenne
sequence}~$M=(M_n)$ defined by
\[
M_n=2^n-1\medspace\medspace (n=1,2,\dots)
\]
contains infinitely many prime terms, which are known as
\emph{Mersenne primes}.

The Mersenne prime conjecture is related to a classical problem in
number theory concerning \emph{perfect numbers}. A whole number is
said to be perfect if, like~$6=1+2+3$ and~$28=1+2+4+7+14$, it is
equal to the sum of all its proper divisors. Euclid pointed out
that~$2^{k-1}(2^k-1)$ is perfect whenever~$2^k-1$ is prime. A much
less obvious result, due to Euler, is a partial converse: if~$n$ is
an \emph{even} perfect number, then it must have the
form~$2^{k-1}(2^k-1)$ for some~$k$ with the property that~$2^k-1$ is
a prime. Whether there are any \emph{odd} perfect numbers remains an
open question. Thus finding Mersenne primes amounts to finding
(even) perfect numbers.

The sequence~$M$ certainly produces some primes initially, for
example,
$$
M_2=3,M_3=7,M_5=31,M_7=127,\dots.
$$
However, the appearance of Mersenne primes quickly thins out: only
forty-three are known, the largest of which,~$M_{30,402,457}$, has
over nine million decimal digits. This was discovered by a team at
Central Missouri State University as part of the GIMPS
project~\cite{gimps}, which harnesses idle time on thousands of
computers all over the world to run a distributed version of the
Lucas--Lehmer test.

A paltry forty-three primes might seem rather a small return for
such a huge effort. Anybody looking for gold or gems with the same
level of success would surely abandon the search. It seems fair to
ask why we should expect there to be infinitely many Mersenne
primes. In the absence of a rigorous proof, our expectations may be
informed by \emph{heuristic} arguments. In section~3 we discuss
heuristic arguments for this and other more or less tractable
problems in number theory.

\medskip
\noindent{\bf Primitive prime divisors.} In 1892, Zsigmondy~\cite{Z}
discovered a beautiful argument that shows that the sequence~$M$
does yield infinitely many prime numbers---but in a less restrictive
sense. Given any integer sequence~$S=(S_n)_{n\geq 1}$, we define a
\emph{primitive divisor} of the term~$S_n$ $(\neq 0)$ to be a
divisor of~$S_n$ that is coprime to every nonzero term~$S_m$
with~$m<n$. Any prime factor of a primitive divisor is called a
\emph{primitive prime divisor}. Factoring the first few terms of the
Mersenne sequence reveals several primitive divisors, shown in bold
in Table~\ref{tableofprimitivedivisors}.
\begin{table}[ht]
\begin{center}
\caption{\label{tableofprimitivedivisors}Primitive divisors
of~$(M_n)$.}
\begin{tabular}{|r|r|c|}
\hline
$n$&$M_n$&Factorization\\
\hline
$2$&$3$&$\boldsymbol{3}$\\
$3$&$7$&$\boldsymbol{7}$\\
$4$&$15$&$3\cdot\boldsymbol{5}$\\
$5$&$31$&$\boldsymbol{31}$\\
$6$&$63$&$3^2\cdot7$\\
$7$&$127$&$\boldsymbol{127}$\\
$8$&$255$&$3\cdot5\cdot\boldsymbol{17}$\\
$9$&$511$&$7\cdot\boldsymbol{73}$\\
$10$&$1023$&$3\cdot\boldsymbol{11}\cdot31$\\
\hline
\end{tabular}
\end{center}
\end{table}
Notice that the term~$M_6$ has no primitive divisor, but each of the
other early terms has at least one. Zsigmondy~\cite{Z} proved that
all the terms~$M_n$~$(n>6)$ have primitive divisors. He also proved
a similar result for more general sequences~$U=(U_n)_{n\ge 1}$,
namely, those of the form~$U_n=a^n-b^n$, where~$a$ and~$b$~$(a>b)$
are positive coprime integers:~$U_n$ has a primitive divisor
unless~$a=2,b=1$ and~$n=6$ or~$a+b$ is a power of~$2$ and~$n=2.$

Apart from the special situation in which~$a-b=1$, it is not
reasonable to expect the terms~$U_n=a^n-b^n$ ever to be prime, since
the identity
\[
{a^n-b^n}=(a-b)(a^{n-1}+a^{n-2}b+a^{n-3}b^2+\cdots+b^{n-1})
\]
shows that~$U_n$ is divisible by~$a-b$. However, it does seem
likely that for any coprime starting values~$a$ and~$b$ infinitely
many terms of the sequence~$\left(U_n/(a-b)\right)$ might be prime.
Sadly, no proof of this plausible statement is known for even a
single pair of starting values.

Although Zsigmondy's result is much weaker than the Mersenne prime
conjecture, it initiated a great deal of interest in the arithmetic
of such sequences (see~\cite[chap.~6]{MR2004c:11015}). It has also
been applied in finite group theory (see
Praeger~\cite{MR2000h:20090}, for example).
Schinzel~\cite{MR26:1280},~\cite{MR49:8961} extended Zsigmondy's
result, giving further insight into the finer arithmetic of
sequences like~$M$. For example, he proved that~$M_{4k}$ has a
composite primitive divisor for all odd~$k$ greater than five.

\bigskip

\addtocounter{section}{1}
\noindent{\bf\thesection\ RECURRENCE
SEQUENCES.}\setcounter{theorem}{0}
For most people their
first introduction to the Fibonacci sequence
$$
A_1=1,A_2=1,A_3=2,A_4=3,A_5=5,A_6=8,\dots
$$
is through the (binary) linear recurrence relation
$$
A_{n+2}=A_{n+1}+A_n.
$$
Sequences such as the Mersenne sequence~$M$ and those considered
by Zsigmondy are of particular interest because they also satisfy
binary recurrence relations. The terms~$U_n=a^n-b^n$ satisfy the
recurrence
$$
U_{n+2}=(a+b)U_{n+1}-abU_n\medspace\medspace(n=1,2,\dots).
$$

More generally, let~$u$ and~$v$ denote conjugate quadratic integers
(i.e., zeros of a monic irreducible quadratic polynomial with
integer coefficients). Consider the integer sequences~$U(u,v)$
and~$V(u,v)$ defined by
$$
U_n(u,v)=(u^n-v^n)/(u-v),\medspace\medspace V_n(u,v)=u^n+v^n.
$$
For instance, the Fibonacci sequence is given by
$$
A_n=U_n\left(\frac{1+\sqrt{5}}{2}, \frac{1-\sqrt{5}}{2}\right)\!.
$$
The sequence~$U(u,v)$ satisfies the recurrence relation
$$
U_{n+2}=(u+v)U_{n+1}-uvU_n\medspace\medspace(n=1,2,\dots),
$$
and~$V(u,v)$ satisfies the same relation.

Some powerful generalizations of Zsigmondy's theorem have been
obtained for these sequences. Bilu, Hanrot, and
Voutier~\cite{MR2002j:11027} used methods from Diophantine analysis
to prove that both~$U_n(u,v)$ and~$V_n(u,v)$ have primitive divisors
once~$n>30$. The two striking aspects of this result are the uniform
nature of the bound and its small numerical value. In particular,
for any given sequence it is easy to check the first thirty terms
for primitive divisors, arriving at a complete picture. For example,
an easy calculation reveals that the Fibonacci number~$A_n$ does not
have a primitive divisor if and only if~$n=1,2,6,$ or~$12$.

\medskip

\noindent{\bf Bilinear recurrence sequences.} The theory of linear
recurrence sequences has a \emph{bilinear} analogue. For example,
the Somos-4 sequence~$S=(S_n)$ is given by the bilinear recurrence
relation
$$
S_{n+4}^{\vphantom{2}}S_n^{\vphantom{2}} =
S_{n+3}^{\vphantom{2}}S_{n+1}^{\vphantom{2}}+S_{n+2}^2\medspace\medspace(n=1,2,\dots),
$$
with the initial condition~$S_1=S_2=S_3=S_4=1$.
This sequence begins
$$
1,1,1,1,2,3,7,23,59,314,1\thinspace529,
8\thinspace209,833\thinspace313,620\thinspace297,
7\thinspace869\thinspace898,\dots.
$$
Amazingly, all the terms are integers even though
calculating~$S_{n+4}$ {a priori} involves dividing by~$S_n$. This
sequence was discovered by Michael Somos~\cite{somoscrux}, and it is
known to be associated with the arithmetic of elliptic curves
(see~\cite[secs.~10.1,~11.1]{MR2004c:11015} for a summary of this,
and further references, including a remarkable observation due to
Propp {et al.} that the terms of the sequence must be integers
because they count matchings in a sequence of graphs.)

Amongst the early terms of~$S$ are several primes: of those that we
listed,
$$
2,3,7,23,59,8\thinspace209,620\thinspace297
$$
are prime. It seems natural to ask whether there are infinitely many
prime terms in the Somos-4 sequence. More generally, consider
integer sequences~$S$ satisfying relations of the type
\begin{equation}\label{equation:generalsomos}
S_{n+4}S_n=eS_{n+3}S_{n+1}+fS_{n+2}^2,
\end{equation}
where~$e$ and~$f$ are integral constants not both zero. Such
sequences are often called \emph{Somos sequences} (or \emph{bilinear
recurrence sequences}) and Christine Swart~\cite{swart}, building on
earlier remarks of Nelson Stephens, showed how they are related to
the arithmetic of elliptic curves. Some care is needed because, for
example, a binary linear recurrence sequence always satisfies some
bilinear recurrence relation of this kind. We refer to a Somos
sequence as \emph{nonlinear} if it does not satisfy any linear
recurrence relation. These are natural generalizations of linear
recurrence sequences, so perhaps we should expect them to contain
infinitely many prime terms. Computational evidence
in~\cite{MR88h:11094} tended to support that belief because of the
relatively large primes discovered. However, a heuristic argument
(discussed later) using the prime number theorem was adapted
in~\cite{EEWelliptic}, and it suggested that a nonlinear Somos
sequence should contain only finitely many prime terms.
See~\cite{MR2045409} for proofs in some special cases.

On the other hand, Silverman~\cite{MR89m:11027} established a
qualitative analogue of Zsigmondy's result for elliptic curves that
applies, in particular, to the Somos-4 sequence. An explicit form of
this result proved by Everest, McLaren, and Ward~\cite{emw}
guarantees that from~$S_5$ onwards all terms have primitive
divisors. There are many nonlinear Somos sequences to which
Silverman's proof does not apply. A version of Zsigmondy's theorem
valid for these sequences awaits discovery.

\medskip

\noindent{\bf Polynomials.} Given the previous sections, it might be
tempting to think that all integral recurrence sequences have
primitive divisors from some point on. However, it is easy to write
down counterexamples. The sequence~$T=(T_n)$ defined by~$T_n=n$,
which satisfies
$$
T_{n+2}=2T_{n+1}-T_n,
$$
is a binary linear recurrence sequence that does not always produce
primitive divisors. This is a rather trivial counterexample, so
consider now the sequence~$P$ defined by
$$
P_n=n^2 + \beta,
$$
where~$\beta$ is a nonzero integer. The terms of this sequence
satisfy the linear recurrence relation
$$
P_{n+3}=3P_{n+2}-3P_{n+1}+P_n. 
$$
It has long been suspected that for any fixed~$\beta$ such
that~$-\beta$ is not a square the sequence~$P$ contains infinitely
many prime terms. A proof is known for not even one value
of~$\beta$. It seems reasonable to ask the apparently simpler
question about the existence of primitive divisors of terms. Clearly
any prime term is itself a primitive divisor, but do the composite
terms have primitive divisors? Using a result of Schinzel about the
largest prime factor of the terms in polynomial sequences it is
fairly easy to prove the following:

\begin{theorem}\label{maincorollary}
If~$-\beta$ is not a square, then there are infinitely many terms of
the sequence~$P$ that do not have primitive divisors.
\end{theorem}

We prove Theorem~\ref{maincorollary} in section~4. Computations
suggest that the following stronger result should be true.

\begin{conjecture}\label{conjecture}
Suppose that~$-\beta$ is not a square. If~$\rho_{\beta}(N)$ denotes
the number of terms~$P_n$ in the sequence~$P$ with~$n<N$ that have
primitive divisors, then
$$
\rho_{\beta}(N)\sim cN
$$
for some constant~$c$ satisfying~$0<c<1$.
\end{conjecture}
Here, and throughout the article, for functions~$f:\mathbb
R\rightarrow\mathbb R$ and~$g:\mathbb R\rightarrow\mathbb R_+,$ we
write~$f\sim g$ to mean~$f(x)/g(x)\rightarrow 1$ as~$x\rightarrow
\infty$.

In the last section of the article, we consider some approaches to
bounding the number of terms in~$P$ that have primitive divisors.
For example, we will furnish a simple proof that
$$
\liminf_{N\to\infty}\frac{\rho_{\beta}(N)}{N}\ge\frac{1}{2}.
$$
We have been unable to find a proof of Conjecture~\ref{conjecture}.
In section~$3$ we show how other kinds of arguments can be
marshalled in its support, and in section~$4$ we discuss briefly the
nature of the constant~$c$.

\medskip

\noindent{\bf Linear recurrence sequences.} To set matters in a more
general context, define~$L=(L_n)_{n\ge 1}$ to be a \emph{linear
recurrence sequence of order~$k$} ($k\ge1$) if it satisfies a
relation
\begin{equation}\label{genrecrel}
L_{n+k}=c_{k-1}L_{n+k-1}+\dots
+c_0L_n\medspace\medspace(n=1,2,\dots)
\end{equation}
for constants~$c_0,\dots,c_{k-1}$, but satisfies no shorter
relation. When~$k=3$ (respectively,~$k=4$), the sequence~$L$ is
called a \emph{ternary (respectively, quaternary) linear recurrence
sequence}. For example, the sequences~$P$ considered in the previous
section are all ternary linear recurrence sequences.
Theorem~\ref{maincorollary} shows that Zsigmondy's theorem cannot
extend to these quadratic sequences. Some nonpolynomial sequences
that cannot satisfy Zsigmondy will now be presented.

With~$u$ and~$v$ again denoting conjugate quadratic integers, the
integer sequence~$W(u,v)=(W_n(u,v))_{n\ge 1}$ defined by
$$
W_n(u,v)=(u^n-1)(v^n-1)
$$
is always a linear recurrence sequence.

\begin{example}\label{BandC}
The sequence~$B=-W(2+\sqrt{3},2-\sqrt{3})$ begins
$$2,12,50,192,722,2700,10082,37632,140450,524172,\dots,
$$
and it is a ternary sequence satisfying
$$B_{n+3}=5B_{n+2}-5B_{n+1}+B_n.
$$
From the seventh term on, all the terms of the
sequence seem to have primitive divisors.
\end{example}

\begin{example}\label{BandC2}
The sequence~$C=-W(1+\sqrt{2},1-\sqrt{2})$ begins
$$
2,4,14,32,82,196,478,1152,2786,6724,\dots,
$$
and it is a quaternary sequence satisfying
$$
C_{n+4}=2C_{n+3}+2C_{n+2}-2C_{n+1}-C_n.
$$
In contrast to the previous example, the terms~$C_{2k}$ for odd~$k$
do not have primitive divisors.
\end{example}

In general, when~$uv=-1$, the terms~$W_{2k}(u,v)$ for odd~$k$ fail
to yield primitive divisors. This is because an easy calculation
reveals that
$$
W_{2k}(u,v)=-W_k(u,v)^2
$$
when~$k$ is odd. On the other hand, we recommend the following as an
exercise: when~$uv=1$, the terms~$W_n(u,v)$ do produce primitive
divisors from some point on. As far as we can tell, to establish
this requires Schinzel's extension~\cite{MR49:8961} of Zsigmondy's
result to the algebraic setting. (We are indebted to Professor
Gy\"ory for communicating to us the remarks about~$W(u,v)$.)

All of these special cases can be subsumed into a wider picture. Write
$$
f(x)=x^k-c_{k-1}x^{k-1}-\dots -c_0
$$
for the \emph{characteristic polynomial} of the linear recurrence
relation in~\eqref{genrecrel}. Then~$f$ can be factored
over~$\mathbb C$,
$$
f(x)=(x-\alpha_1)^{e_1}\dots (x-\alpha_d)^{e_d}.
$$
The algebraic numbers~$\alpha_1,\dots ,\alpha_d$, are known as the
\emph{characteristic roots} (or just \emph{roots}) of the
sequence. The terms~$L_n$ of any sequence~$L$ satisfying the
relation~\eqref{genrecrel} can be written
$$
L_n=\sum_{i=1}^dg_j(n)\alpha_i^n
$$
for polynomials~$g_1,\dots,g_d$ of degrees~$e_1-1,\dots ,e_d-1$
with algebraic coefficients.

The roots of the sequences in Examples~\ref{BandC} and~\ref{BandC2}
are quite different in character. In general, if~$uv=1$,
then~$W(u,v)$ is a ternary linear recurrence sequence with
roots~$1,u,$ and~$v$. When~$uv=-1$,~$W(u,v)$ is quaternary with
roots~$1,-1,u,$ and~$v$. For the quadratic sequence defined
by~$P_n=n^2+\beta$,~$\alpha=1$ is a triple root of the associated
characteristic polynomial.

It seems reasonable to conjecture that the terms of an integral
linear recurrence sequence of order greater than one will have
primitive divisors from some point on provided that its roots are
distinct and no quotient~$\alpha_i/\alpha_j$ of different roots is a
root of unity.

\bigskip

\addtocounter{section}{1}
\noindent{\bf\thesection\ HEURISTIC
ARGUMENTS.}\setcounter{theorem}{0}\label{heuristicssection}
There are a number of ways that mathematicians form a view on which
statements are likely to be true. These views inform research
directions and help to concentrate effort on the most fruitful areas
of enquiry.

The only certainty in mathematics comes from rigorous proofs that
adhere to the rules of logic: the discourse of {\it logos}. When
such a proof is not available, other kinds of arguments can make
mathematicians expect that statements will be true, even though
these arguments fall well short of proofs. These are called
\emph{heuristic} arguments---the word comes from the Greek root
\foreignlanguage{greek}{\it Eurhka} ({\it Eureka}), meaning ``I have
found it.'' It usually means the principles used to make decisions
in the absence of complete information or the ability to examine all
possibilities. In informal ways, mathematicians use heuristic
arguments all the time when they discuss mathematics, and these are
part of the {\it mythos} discourse in mathematics.

One consequence of the prime number theorem is the following
statement: the probability that~$N$ is prime is roughly~$1/\log N$.
What this means is that if an element of the set~$\{1,\dots,N\}$ is
chosen at random using a fair~$N$-sided die, then the
probability~$\rho_N$ that the number chosen is prime
satisfies~$\rho_N\log N\to 1$ as~$N\to\infty$. This crude estimate
has been used several times to argue heuristically in favor of the
plausibility of conjectured solutions of difficult problems. Some
examples follow. In each case the argument presented falls well
short of a proof, yet it still seems to have some predictive power
and has suggested lines of attack.

\medskip

\noindent{\bf Fermat primes.} Hardy and
Wright~\cite[sec.~2.5]{MR568909} argued along these lines that there
ought to be only finitely many Fermat primes. A \emph{Fermat prime} is a
prime number in the sequence~$(F_n)$ of \emph{Fermat numbers}:
$$
F_n=2^{2^n}+1.
$$
Fermat demonstrated that~$F_1,F_2,F_3,$ and~$F_4$ are all primes;
Euler showed that~$F_5$ is composite by using congruence arguments.
Since then, many Fermat numbers have been shown to be composite and
quite a few have been completely factored. Wilfrid Keller maintains
a web site~\cite{fermatsite} with details of the current state of
knowledge on factorization of Fermat numbers. The number of Fermat
primes~$F_n$ with~$n<N$, if they are no more or less likely to be
prime than a random number of comparable size, should be roughly
$$
\sum_{n<N}\frac{1}{\log F_n}\sim\sum_{n<N}\frac{1}{2^n\log
2}<\frac{1}{\log2}.
$$
Statements like this cannot be taken too literally, for the
numbers~$F_n$ have many special properties, not all of which are
understood. However, this kind of argument tends to support the
belief that there are only finitely many Fermat primes and would
incline many mathematicians to attempt to prove that statement
rather than its negation. Massive advances in computing power
suggest that we know---indeed, that Fermat knew---all the Fermat
primes.

\medskip

\noindent{\bf Mersenne primes.} The prime number theorem can also be
used to argue in support of the Mersenne prime conjecture. A
heuristic argument of the following form is used. First,~$2^k-1$ can
be prime only for~$k$ a prime, so assume now that~$k$ is a
prime~$p$. We would like to estimate the probability that~$2^p-1$ is
prime. The prime number theorem suggests that a random number of the
size of~$2^p-1$ is prime with probability~$1/\log(2^p-1)$, which is
around~$1/p\log2$. However,~$2^p-1$ is far from random: it is not
divisible by~$2$, nor by~$3$, and indeed not by any prime smaller
than~$2p$. Arguing in this way suggests that the
probability that~$2^p-1$ is prime is approximately
\begin{equation}\label{equation:mersenneheuristic}
\rho_p=\frac{1}{p\log 2}
\cdot\frac{2}{1}\cdot\frac{3}{2}\cdot\frac{5}{4}\cdot\cdots
\cdot\frac{q}{q-1},
\end{equation}
where~$q$ is the largest prime less than~$2p$. This suggests that
the expected number of Mersenne primes~$M_n$ with~$n<N$ is
roughly~$\sum_{p<N}\rho_p.$

Since~$\rho_p>1/p\log2$, the sum diverges by Mertens's theorem
(see~\eqref{equationmertensformula} for a precise statement), which
suggests that there are infinitely many Mersenne primes.
Wagstaff~\cite{MR84j:10052} and then Pomerance and
Lenstra~\cite{MR873581} have extended this heuristic argument by
including estimates for the product of rationals
in~\eqref{equation:mersenneheuristic} to obtain an asymptotic
estimate that closely matches the available evidence. On the basis
of these heuristics, they conjecture that the number of Mersenne
primes~$M_n$ with~$n<N$ is asymptotically
$$
\frac{e^{\gamma}}{\log 2}\log N,
$$
where~$\gamma$ is the Euler-Mascheroni constant. Caldwell's Prime
Page~\cite{primepages} gives
more details about these arguments and about the hunt for new Mersenne
primes.

\medskip

\noindent{\bf Bilinear recurrence sequences.} Consider now the Somos
sequences defined by the recurrence~\eqref{equation:generalsomos}.
General results about heights on elliptic curves show that the
growth rate of~$S_n$ is quadratic-exponential. In other words,
$$
\log S_n \sim hn^2,
$$
where~$h$ is a positive constant. Thus the expected number of
prime terms with~$n<N$ should be approximately
$$
\sum_{n<N}\frac{1}{\log
S_n}\sim\frac{1}{h}\sum_{n<N}\frac{1}{n^2}\le\frac{\pi^2}{6h}.
$$
This resembles the argument of Hardy and Wright and suggests that
only finitely many prime terms should be expected. Proofs of the
finiteness in many special cases have subsequently been
found~\cite{MR2045409}. The search for these proofs was motivated in
part by the heuristic arguments. Interestingly, it is known that the
constant~$h$ is uniformly bounded below across all nonlinear
integral Somos sequences. Thus the style of this heuristic argument
suggests that perhaps the total number of prime terms is uniformly
bounded across all such sequences. Extensive calculation has failed
to yield a sequence with more than a dozen prime terms.

\medskip

\noindent{\bf Quadratic polynomials.} Suppose~$\beta$ is an integer
that is not the negative of a square and recall the sequence~$P$
given by~$P_n=n^2+\beta$. Again, the prime number theorem predicts
that there are roughly
$$
\sum_{n<N}\frac{1}{\log P_n}
$$
prime terms in the sequence~$P$ with~$n<N$, assuming again
that~$P_n$ is neither more nor less likely to be prime than a random
number of that size. The sum is asymptotically~$N/(2\log N)$, which
supports the belief that there are infinitely many prime terms in
the sequence~$P$. Computation suggests that for fixed~$\beta$ there
will be~$dN/\log N$ prime terms with~$n<N$, where~$d=d(\beta)$ is a
constant that depends upon~$\beta$. Bateman and
Horn~\cite{MR0148632} offered a heuristic argument and provided
numerical evidence to suggest that
$$
d=\frac{1}{2}\,{\displaystyle\prod_{p}}\left(\vphantom{\sum}\right.1-
\frac{1}{p}\left.\vphantom{\sum}\right)^{-1}
\left(\vphantom{\sum}\right.1-
\frac{w(p)}{p}\left.\vphantom{\sum}\right),
$$
where the product is taken over all primes and~$w(p)$ denotes the
number of solutions~$x$ modulo~$p$ to the congruence~$x^2\equiv
-\beta\pmod{p}$.

\bigskip

\addtocounter{section}{1}
\noindent{\bf\thesection\ BIASED
NUMBERS.}\setcounter{theorem}{0}\label{proofs}
We now return to the problem of looking for primitive prime factors
in the sequence given by~$P_n=n^2+\beta$ with~$-\beta$ not a square.
Since we are mainly interested in asymptotic behaviour, we assume
from now on that~$n>\vert\beta\vert$. The terms~$P_n$
with~$n\le\vert\beta\vert$ are not guaranteed to exhibit the
behavior described in this section.

\begin{lemma}\label{lemmaaddedbyreferee}
A prime~$p$ is a primitive divisor of~$P_n$ if and only if~$p$
divides~$P_n$ and~$p>2n$.
\end{lemma}

\begin{proof}
Consider first a prime~$p$ dividing~$P_n$ with~$p< n$. Then, by
assumption,~$P_n\equiv 0\pmod{p}$, so~$P_m\equiv 0\pmod{p}$ for
some~$m$ smaller than~$p$ simply by choosing~$m$ to be the residue
of~$n$ modulo~$p$. Because~$p< n$,~$m<n$. In others words,~$p$ is
not a primitive divisor of~$P_n$.

This means that to find primitive divisors of~$P_n$ we have to look
for prime divisors that are greater than~$n$. (Note that~$n$ does not
divide~$P_n$, as~$n>|\beta|$.) We can say more: we
can guarantee a solution of~$P_m\equiv 0\pmod{p}$ for some~$m$
satisfying~$m\le p/2$. Thus, to find primitive divisors we have to
look only amongst the prime divisors that are bigger than~$2n$ (i.e.
a prime~$p$ dividing~$P_n$ is a primitive divisor only if~$p>2n$).

Conversely, suppose that~$p$ is a prime dividing~$P_n$ that is not a
primitive divisor. Then~$n^2+\beta \equiv0\pmod{p}$, and there is an
integer~$m$ ($<n$) with~$m^2+\beta\equiv0\pmod{p},$ so (by
subtracting the two congruences)~$m^2-n^2\equiv0\pmod{p}$. It
follows that~$m\pm n\equiv0\pmod{p}$. In particular,
$$
p\le m+n<2n.
$$
It follows that a prime~$p$ is a primitive divisor of~$P_n$ if and
only if~$p$ divides~$P_n$ and~$p>2n$.
\end{proof}

We call an integer~$k$ is \emph{\primeish} if it has a prime factor~$q$
with~$q>2\sqrt k$. Thus any prime greater than three is \primeish.
The numbers~$22,26,$ and~$34$ are \primeish, whereas~$24$ and~$28$
are not.

\begin{proposition}\label{primeishiff} When~$n>|\beta|$,
the term~$P_n$ has a primitive divisor if and only if~$P_n$ is
\primeish. If~$n>|\beta|$ and~$P_n$ has a primitive divisor, then
that primitive divisor is a prime, and it is unique.
\end{proposition}

\begin{proof}
Part of the first statement comes from
Lemma~\ref{lemmaaddedbyreferee}. To complete the proof of the first
statement we claim that, for~$n$ greater than~$\vert \beta
\vert$,~$P_n$ has such a prime divisor if and only if~$P_n$ is
\primeish. If~$p$ is a prime dividing~$P_n$ and~$p>2n$, then
$$
p\ge2n+1>2\sqrt{n^2+n}>2\sqrt{n^2+\beta}.
$$
Conversely, if~$p>2\sqrt{n^2+\beta}$, then
$$
p\ge2\sqrt{n^2-n+1}>2n-1,
$$
so~$p>2n$ (since~$2n$ cannot be prime).

The uniqueness of the primitive divisor follows at once. If~$p$ is a
prime dividing~$P_n$ and~$p>2\sqrt{P_n}$, then no other prime
divisor can be as large, hence cannot be primitive.
\end{proof}

The requirement~$n>\vert \beta\vert$ is necessary: if~$|\beta|$ is
prime, then~$P_{|\beta|}$ has primitive divisor~$|\beta|$ but is not
\primeish. Also, terms with small~$n$ may have more than one
primitive divisor. For example, the sequence of values of the
polynomial~$n^2+6$ begins~$7,10,\dots$ so the second term has two
primitive prime divisors. The kind of results discussed here are
asymptotic results, which makes this restriction unimportant.

\begin{proof}[Proof of Theorem~\ref{maincorollary}]
Results of Schinzel~\cite[Theorem~13]{MR0222034} show that for any
positive~$\alpha$, the largest prime factor of~$P_n$ is bounded
above by~$n^{\alpha}$ for infinitely many~$n$. Taking $\alpha={1}$,
we conclude that~$P_n$ is not biased infinitely often. By
Proposition~\ref{primeishiff},~$P_n$ fails to have a primitive
divisor infinitely often.
\end{proof}

In section~$5$ of the article, we consider quantitative
information about the frequency with which, rather than the extent
to which,~$P_n$ is not biased.

Support for Conjecture~\ref{conjecture} follows from
Proposition~\ref{primeishiff} because an asymptotic
formula can be obtained
for the distribution of \primeish\ numbers. Alongside the earlier
notation for describing the growth rates of various functions, we
also use the following:
Given functions~$f:\mathbb R\rightarrow\mathbb R$
and~$g:\mathbb R\rightarrow\mathbb R_+,$
we write~$f=O(g)$ to
mean that~$|f(x)|/g(x)$ is bounded and~$f=o(g)$ to signify
that~$f(x)/g(x)\rightarrow 0$ as~$x\rightarrow\infty$.

\begin{theorem}\label{lopcount}
If~$\pi_{\imath}(N)$ denotes the number of \primeish\ numbers less
than or equal to~$N$, then
$$
\pi_{\imath}(N)\sim N\log 2.
$$
\end{theorem}

\begin{proof} Write a
\primeish\ number as~$qm$, where~$q$ is its largest prime factor.
The \primeish\ condition then translates to~$q>4m$. To compute the number
of \primeish\ numbers below~$N$, note that the counting can be
achieved by dividing the set into two parts. Let~$p$ denote a
variable prime. When~$p<2\sqrt N$ there are~$\lfloor p/4\rfloor$
\primeish\ integers~$pm$ smaller than~$N$ (here~$\lfloor x\rfloor$ denotes
the greatest integer less than or equal to~$x$).
When~$p\ge 2\sqrt N$, each number~$pm$ smaller than~$N$ is
\primeish, so there are~$\lfloor{N}/{p}\rfloor$ \primeish\
integers~$pm$. Hence the total number is
$$
\sum_{p<2\sqrt
N}\left\lfloor\textstyle\frac{p\vphantom{N}}{4\vphantom{p}}
\right\rfloor+\negmedspace\negmedspace\sum_{2\sqrt N\le
p<N}\left\lfloor\textstyle\frac{N}{p}\right\rfloor.
$$
The first sum is~$O(N/\log N)$ and can be ignored
asymptotically. The second sum differs from
$$
N\negmedspace\negmedspace\negmedspace\negmedspace\sum_{2\sqrt N\le
p<N}\negmedspace\negmedspace\negmedspace\textstyle\frac{1}{p}
$$
by an amount that is~$O(N/\log N)$ by the prime number theorem.
To estimate this sum we use Mertens's formula, which can be found in
Apostol's book~\cite[Theorem~4.12]{MR0434929}:
\begin{equation}\label{equationmertensformula}
\sum_{p<x}\frac{1}{p}=\log \log x + A + o(1).
\end{equation}
Applying~\eqref{equationmertensformula} thus estimates~$\pi_{\imath}(N)$
as
$$
N\left[\log\log N+A-\log(\log\sqrt N+\log 2)-A+o(1)\right],
$$
which is asymptotically~$N\log 2$.
\end{proof}

Theorem~\ref{lopcount} can be applied to give the following
heuristic argument in support of Conjecture~\ref{conjecture}. The
probability that a large integer is \primeish\ is roughly~$\log
2$. Hence the expected number of \primeish\ values of~$n^2+\beta$
with~$n<N$ is asymptotically~$N\log 2$. Computational evidence
suggests that the number of \primeish\ terms in~$n^2+\beta$ is
asymptotically~$cN$ for some constant~$c$. Computations
with~$|\beta|<10$ suggest the constant~$c$ looks reasonably close
to~$\log 2$ in each case, although convergence appears slow.

\bigskip

\addtocounter{section}{1}
\noindent{\bf\thesection\ COUNTING PRIMITIVE DIVISORS.}\setcounter{theorem}{0}
The article concludes with some simple estimates
for~$\rho_{\beta}(N)$, the number of terms~$P_n$
in the sequence~$P$ with~$n<N$ that have
primitive divisors. The proofs use little aside
from well-known estimates for sums over primes, which can be found
in the book of Apostol~\cite{MR0434929}.

\begin{theorem}\label{infinitelymany}
There is a constant~$C>0$ such that
\begin{equation}\label{firsttheorembound}
\rho_{\beta}(N)<N-\frac{CN}{\log N}
\end{equation}
holds for all sufficiently large~$N$.
There is a constant~$D>0$ such that
\begin{equation}\label{secondtheorembound}
\frac{N}{2}-\frac{DN}{\log N}<\rho_{\beta}(N)
\end{equation}
is true for all sufficiently large~$N$.
\end{theorem}

Both of the statements in Theorem~\ref{infinitelymany} can be
strengthened along the following lines: any choice of
constants~$C$ or~$D$ could be made.
As each constant varies, so does the smallest value of~$N$ beyond
which the inequalities become valid.

Apart from a finite number of primes, any prime~$p$ that
divides~$n^2+\beta$ has the property that~$-\beta$ is a quadratic
residue modulo~$p$. Let~$\cR$ denote the set of odd primes for
which~$-\beta$ is a quadratic residue. Notice that~$\cR$ comprises
the intersection of a finite union of arithmetic progressions with
the set of primes and that this finite union of arithmetic progressions
in turn comprises exactly half of the residue classes
modulo~$4\vert\beta\vert$. We will prove the two parts of
Theorem~\ref{infinitelymany} in reverse order, because the upper
bound~\eqref{firsttheorembound} arises by specializing the argument
used to prove the lower bound~\eqref{secondtheorembound}.

Write
$$
Q_N=\prod_{n=1}^N\vert P_n\vert,
$$
and denote by~$\omega (Q_N)$ the number of distinct prime divisors
of~$Q_N$. By Proposition~\ref{primeishiff} it is sufficient to
bound~$\omega (Q_N)$ because, with finitely many exceptions, a
primitive divisor is unique.

\medskip

\noindent{\bf Proof of the lower bound.}
Define
\[
\cS=\{p\in \cR\mid p\vert Q_N\mbox{ and }p<2N\},\medspace
\cS'=\{p\in
\cR\mid p\vert Q_N\mbox{ and }p\ge 2N\}.
\]
Let~$s=\vert \cS\vert$ and~$s'=\vert \cS'\vert$. We seek a lower
bound for~$s+s'$, since
\[
s+s'=\omega(Q_N).
\]
The asymptotic form of Dirichlet's theorem\footnote{In~$1826$
Dirichlet proved that if~$a$ and~$b$ are positive integers with no
common factor, then there are infinitely many primes of the
form~$ax+b$ with~$x$ in~$\mathbb N$. This result, which appeared in a
memoir published in 1837~\cite{dirichlet}, was proved using methods
from analysis, thus laying the foundations for the subject now
called analytic number theory. Writing~$\pi_a(X)$ for the number of
primes of the form~$ax+b$ with~$x<X$, Dirichlet proved
that~$\pi_a(X)\to\infty$ as~$X\to\infty$. There is also what might
be called a prime number theorem for arithmetic progressions, which
gives an asymptotic estimate for the number of such primes. It
states that~$\pi_a(X)\sim X/{\phi(a)}\log X$, where~$\phi$ is the
Euler totient function. This was shown by de la Vall{\'e}e
Poussin; a proof can be found in the book of
Prachar~\cite[chap.~5, sec.~7]{prachar}. It is this result that we are
using here.} on primes in arithmetic progression implies that
asymptotically half the primes lie in~$\mathcal R$, so
\begin{equation}\label{asymptoticfors} s\sim\frac{N}{\log
2N}.
\end{equation}
Therefore it is sufficient to estimate~$s'$ from below.

\begin{proof}[Proof of equation~\eqref{secondtheorembound}]
From the definition of~$Q_N$,
$$
\log Q_N=\sum_{n=1}^N\log\vert n^2+\beta\vert
=2\sum_{n=1}^{N}\left(\log n+\vphantom{A^{A^A}}
O\negthinspace\left(\textstyle\frac{1}{n^2}\right)\right)
=\left(2\sum_{n=1}^{N}\log n\right)\negthinspace+O(1),
$$
so by Stirling's formula
\begin{equation}\label{estimate1}
\log Q_N= 2N\log N - 2N + O(1).
\end{equation}
On the other hand, we can write
\begin{equation}\label{applstirling}
\sum_{p|Q_N}e_p\log p=\log Q_N,
\end{equation}
for positive integers~$e_p$
corresponding to the prime decomposition~$\prod_{p\vert Q_N}
p^{e_p}$ of~$Q_N$. The first step in
the proof is to identify a subset of~$\cR$ that contributes a fixed
amount to the main term in~\eqref{estimate1}. The sum on the
left-hand side of~\eqref{applstirling} can be decomposed to give
\begin{equation}\label{estimate2}
\sum_{p\in \cS, p<N}e_p\log p + \sum_{p\in\cS,p\ge N}e_p\log p
+\sum_{p\in \cS'}\log p=\log Q_N,
\end{equation}
noting that~$e_p=1$ whenever~$p\ge 2N$. The second term in the
decomposition is~$O(N)$, since~$e_p\le2$ for~$p$ in~$\cS$ with~$p>N$,
each term in the sum is no larger than~$\log 2N$, and the
prime number theorem implies that there are~$O(N/\log N)$ terms.
Thus the second term does not contribute to the asymptotic
behaviour.

Assume for the moment that
\begin{equation}\label{sumoverpinS}
\sum_{p\in\cS,p<N}e_p\log p=N\log N+O(N).
\end{equation}
Combining~\eqref{estimate1}, \eqref{estimate2},
and~\eqref{sumoverpinS} gives
$$
N\log N + O(N) = \sum_{p\in \cS'}\log p < s'\log P_N=s'\log (N^2+\beta).
$$
Thus~\eqref{secondtheorembound} follows at once, subject to the
proof of~\eqref{sumoverpinS}.
\end{proof}

\begin{proof}[Proof of equation~\eqref{sumoverpinS}]
For each~$p$ in~$\cS$
$$
e_p\ge\left\lfloor\textstyle\frac{2N}{p}\right\rfloor\negmedspace.
$$
Hence there is a constant~$c>0$ such that the left-hand side
of~\eqref{sumoverpinS} is bounded below by
$$
\sum_{p\in \cS,p<N}\left\lfloor\textstyle\frac{2N}{p}\right\rfloor
\log p.
$$
By Apostol~\cite[Theorem~7.3]{MR0434929},
\begin{equation}\label{estimateS}
\sum_{p\in\cS,p<N}\left\lfloor\textstyle\frac{2N}{p}\right\rfloor\log
p = N\log N + O(N).
\end{equation}
For~$p$ in~$\cS$ and~$k$ in~$\N$, denote by~$\ord_p(k)$ the
exponent
of the greatest power of~$p$ dividing~$k$ and put
$$
\cB_p(N)=\{n<N\mid\ord_p(P_n)>1\}.
$$
Then
$$
\sum_{p\in \cS,p<N}e_p\log p = \sum_{p\in \cS,p<N}
\left\lfloor\textstyle\frac{2N}{p}\right\rfloor\log p +
\sum_{p\in\cS,p<N}\left(\vphantom{\sum_A}\right.\sum_{n\in\cB_p(N)}
\ord_p(P_n)-1\left.\vphantom{\sum_A}\right)\log p.
$$
We now show that the second term is asymptotically negligible.
For~$p$ in~$\cS$ the number $-\beta$ has two~$p$-adic square roots,
and~$\ord_p(P_n)=r+1$ if and only if the~$p$-adic expansion of~$n$
agrees with one of these square roots up to the term in~$p^r$ and
no further. Hence
\begin{eqnarray*}
\sum_{p\in\cS,p<N}
\left(\vphantom{\sum_A}\right.\sum_{n\in\cB_p(N)}\ord_p(P_n)-1
\left.\vphantom{\sum_A}\right)\log p &\le&\sum_{p\in\cS,p<N}
\left(\vphantom{\sum_A}\right.\sum_{r=1}^{\frac{\log P_N}{\log p}}
r\cdot 2\left\lceil\textstyle\frac N{p^{r+1}}\right\rceil
\left.\vphantom{\sum_A}\right)\log p \\
&<& 2N\sum_{p\in\cS,p<N}\frac{\log p}{(p-1)^2} + 2s\log P_N,
\end{eqnarray*}
which is~$O(N)$ because the sum converges and~$s=O(N/\log N)$
by~\eqref{asymptoticfors}. Putting this together
with~\eqref{estimateS} establishes~\eqref{sumoverpinS}.
\end{proof}

\medskip

\noindent{\bf Proof of the upper bound.}
This proof is similar to that for the lower bound.
However, it relies on a finer partition of the
set~$\cR$. Given integers~$K>2$ and~$N>K$, we split~$\cS'$ into the
sets
\begin{eqnarray*}
\cT=\{p\in \cR\mid p|Q_N, 2N<p<KN\},\medspace
\cU=\{p\in \cR\mid p|Q_N, KN<p\}.
\end{eqnarray*}

\begin{proof}[Proof of equation~\eqref{firsttheorembound}]
Write~$t=|\cT|$ and~$u=|\cU|$. As before, the contribution
from~$s$ is negligible. Thus we wish to bound the
expression~$t+u$ from above. The sum on the left-hand side
of~\eqref{applstirling} decomposes according to the definitions
of~$\cS,\cT$, and~$\cU$ to give
\[
\sum_{p\in \cS}e_p\log p + \sum_{p\in \cT}\log p
+\sum_{p\in \cU}\log p=\log Q_N,
\]
noting as earlier that~$e_p=1$ whenever~$p>2N$.
Equations~\eqref{estimate1},~\eqref{estimate2},
and~\eqref{sumoverpinS} reveal that
$$
\sum_{p\in \cT}\log p
+\sum_{p\in \cU}\log p<N\log N + aN
$$
for some positive~$a$.
The left-hand side is greater than
$$
t\log N + u\log (KN),
$$
so we add~$t\log K$ to both sides to obtain
$$
(t+u)\log (KN)<N\log N + aN + t\log K.
$$
Rearranging the right-hand side gives
$$
(t+u)\log (KN) < N\log (KN) + (a-\log K)N+t\log K.
$$
Assume that~$K$ is fixed, with~$C=\log K-a>0$.
Dividing through by~$\log(KN)$ leads to
\begin{equation}\label{begindescent}
(t+u)< N - \frac{CN}{\log (KN)} + \frac{t\log K}{\log (KN)}.
\end{equation}
The inequality~$-1/(1+x)<-1+x$ holds when~$x>0$.
We apply this with~$x=\log K/\log N$
to the second term on the right of~\eqref{begindescent}
to obtain the inequality
$$
-\frac{CN}{\log (KN)}<- \frac{CN}{\log (N)}+
O\left(\frac{N}{(\log N)^2}\right),
$$
whose last term is asymptotically negligible.
The last term on the right of~\eqref{begindescent} can be
estimated by appealing to Dirichlet's theorem again,
yielding
$$
\frac{t\log K}{\log (KN)}=O\left(\frac{t}{\log N}\right)=
O\left(\frac{N}{(\log N)^2}\right),
$$
which is also asymptotically negligible. Hence
for any~$C'$ ($0<C'<C$),
$$
\omega(Q_N)\sim t+u<N-\frac{C'N}{\log N}
$$
for all large~$N$.
\end{proof}

A slightly stronger result is provable with these methods, namely,
that
$$
\rho_{\beta}(N)<N-\frac{N\log \log N}{\log N}
$$
for all sufficiently large~$N$.
We leave this as an exercise to the interested reader.


\noindent{\bf GRAHAM EVEREST} received a Ph.D. from King's
College London in 1983 and joined the University of East
Anglia mathematics department the same year. He works in
number theory, specializing in the relationship between
Diophantine problems and elliptic curves. He was ordained
as a minister in the Church of England in 2005.\\
\emph{School of Mathematics, University of East Anglia,
Norwich NR4 7TJ, U.K.\\
g.everest@uea.ac.uk}
\medskip

\noindent{\bf SHAUN STEVENS} received a Ph.D. from King's College
London in 1998. After research positions in Orsay, M{\"u}nster, and
Oxford he joined the University of East Anglia mathematics
department in 2002. He works in representation theory, on the local
Langlands program.\\
\emph{School of Mathematics, University of East
Anglia,
Norwich NR4 7TJ, U.K.\\
shaun.stevens@uea.ac.uk}
\medskip

\noindent{\bf DUNCAN TAMSETT} received a Ph.D.
in Marine Geophysics from the University of Newcastle upon Tyne in 1984.
He is best described as a professional odd-jobs man; currently contributing
to a software system for planning oil well paths, as well as looking after a
system for enhancing side-scan sonar images of the sea-bed and
characterizing/classifying sonar image texture.
He recently moved to
Inverness for the Scottish Highland and Island experience.\\
\emph{5 Drummond Crescent, Inverness IV2 4QW, U.K.\\
duncan@tamsetts.freeserve.co.uk}
\medskip

\noindent{\bf TOM WARD} received a Ph.D. from the University of
Warwick in 1989. After research positions at the University of
Maryland College Park and Ohio State University, he joined the
University of East Anglia mathematics department in 1992. He has
been head of department since 2002. He works in ergodic theory and
its interactions with number theory.\\
\emph{School of Mathematics,
University of East Anglia,
Norwich NR4 7TJ, U.K.\\
t.ward@uea.ac.uk}

\end{document}